\documentclass[10pt,twocolumn,twoside, letterpaper]{IEEEtran} 

\IEEEoverridecommandlockouts 
\usepackage{amsmath,amssymb,graphicx,amsfonts,xspace}
\usepackage{psfrag,enumerate,pdfsync,caption}
\usepackage[all,2cell,dvips]{xy} 
\newcommand{\comments}[1]{}

\usepackage{tikz}
\usetikzlibrary{arrows}
\usetikzlibrary{positioning}
\usetikzlibrary{fit}

\usepackage{xcolor,framed}
\colorlet{shadecolor}{yellow}

\DeclareMathOperator{\diag}{diag}

\DeclareMathOperator{\Span}{span}
\DeclareMathOperator{\image}{Im}
\DeclareMathOperator{\Ker}{Ker}

\newcommand{\SO}{\mathsf{SO}}

\newcommand{\bb}{{\boldsymbol \beta}}
\newcommand{\one}{{\boldsymbol 1}}

\def\qed{\hfill {$\square$}}

\newcommand{\cB}{{\cal B}} 
\newcommand{\cI}{{\cal I}}
 
\newcommand{\cG}{{\cal G}}
\newcommand{\cL}{{\cal L}}

\renewcommand{\Re}{\mathbb{R}} 
\newcommand{\cX}{{\cal X}}

\newcommand{\sR}{{\sf R}} 
\newcommand{\sX}{{\sf X}}
 
\newcommand{\cN}{{\cal N}}
\newcommand{\p}{{A_i}}

\newcommand{\trans}{{}^\top}

\newcommand{\wdes}{\mathbf{f}}
\newcommand{\gdes}{\mathbf{g}}
\newcommand{\hdes}{\mathbf{h}}
\newcommand{\kdes}{\mathbf{r}}
\newcommand{\kq}{k_1}

\newcommand{\prj}{\theta}
\newcommand{\Chi}{X}

\newcommand{\rot}{{\sf rot}}
\newcommand{\tran}{{\sf tran}}

\newtheorem{theorem}{Theorem}
\newtheorem{lemma}{Lemma}

\newtheorem{corollary}{Corollary}
\newtheorem{proposition}{Proposition}
\newtheorem{definition}{Definition}

\newtheorem{claim}{Claim}

\long\def\symbolfootnote[#1]#2{\begingroup%
\def\thefootnote{\fnsymbol{footnote}}\footnote[#1]{#2}\endgroup}

\title{A Smooth Distributed Feedback for Global Rendezvous of Unicycles} 
\author{Ashton Roza, Manfredi Maggiore, Luca Scardovi
\thanks{This research was
  supported by the National Sciences and Engineering Research Council
  of Canada.}
\thanks{The authors are with the Department of Electrical and
  Computer Engineering, University of Toronto, 10 King's College Road,
  Toronto, ON, M5S 3G4, Canada.  {\tt\footnotesize
    ashton.roza@mail.utoronto.ca, maggiore@ece.utoronto.ca, scardovi@scg.utoronto.ca}}}
\begin{document}
\maketitle

\begin{abstract}
This paper presents a solution to the rendezvous control problem for a
network of kinematic unicycles in the plane, each equipped with an onboard
camera measuring its relative displacement with respect to its
neighbors in body frame coordinates.  A smooth, time-independent
control law is presented that drives the unicycles to a common
position from arbitrary initial conditions, under the assumption that
the sensing digraph contains a reverse-directed spanning tree. The
proposed feedback is very simple, and relies only on the onboard
measurements.  No global positioning system is required, nor any
information about the unicycles' orientations.
\end{abstract}

\section{Introduction}\label{sec:intro}

This paper investigates rendezvous control of kinematic
unicycles.  The objective is to design smooth feedbacks for each robot
so as to drive the group to a common position from arbitrary initial
conditions. An important requisite is that the feedback be local and
distributed. In other words, it is required that the feedback depend
only on the relative displacement of each robot to its neighbours
measured in the robot's own body frame, so that the feedback can be
computed using onboard sensing devices such as cameras or laser
systems.

The solution to the rendezvous control problem proposed in this paper
is time-independent and it does not require any information about the
orientation of the unicycles, not even their relative orientation.  To
the best of our knowledge, this is the first solution having the
property of being local and distributed, continuously differentiable,
and time-independent. As we argue below, previous solutions require
either time-varying or discontinuous feedback.  For simplicity of
exposition, the proposed solution relies on the assumption that the
sensing digraph of the unicycles is time-invariant. However, it is
only required to contain a reverse directed spanning tree, which is
the minimal connectivity requirement.

The main difficulty in solving the rendezvous control problem comes from the fact that the unicycles are nonholonomic, in that their velocity is restricted to be parallel to the vehicle's heading direction. To overcome this difficulty, the solution we present relies on a control structure made of two nested loops. An outer loop treats the vehicles as fully-actuated single integrators with a linear consensus controller providing a reference velocity. Here we leverage existing consensus algorithms for single integrators~\cite{beard,moreau,olfati}. The desired velocity computed by the outer loop becomes a reference signal for the inner loop, which assigns local and distributed feedbacks that solve the rendezvous control problem. This methodology is inspired by our previous work in~\cite{RozMag14,roza/online} for rendezvous of rigid bodies in three dimensions.

The rendezvous control problem for unicycles has been investigated
before.  In~\cite{Lin}, the authors presented the first solution.  The
feedback in~\cite{Lin} is local and distributed, but it requires the
use of time-varying feedbacks. In~\cite{dimarogonas} both positions
and attitudes of the unicycles are synchronized using a time-invariant
distributed control. The graph is time-dependent and the authors
assume an initially connected communication graph. The controller that
is implemented, however, is discontinuous. In~\cite{zheng2013} a
time-independent, local and distributed controller is
presented. However, the authors make the assumption that whenever two
vehicles get sufficiently close together they merge into a single
vehicle, introducing a discontinuity in the control function.  To the
best of our knowledge, the solution presented in this paper is the
first one involving feedbacks that are local and distributed,
time-independent, and continuously differentiable. The proposed
solution is of simple implementation, not even requiring any knowledge
about the relative orientation of the unicycles. As we illustrate through simulations, the proposed time-independent, continuously differentiable feedback has practical advantages over the time-varying feedback in~\cite{Lin} and the discontinuous feedback in~\cite{dimarogonas} in that it induces a more natural behaviour in the ensemble of unicycles. The feedback in~\cite{Lin} makes the unicycle ``wiggle'' indefinitely, a behaviour which would be unacceptable in practice. The feedback in~\cite{dimarogonas} induces instantaneous changes in direction that are impossible to achieve with realistic implementations. 


The paper is organized as follows. In Section~\ref{sec:prelims} we
present the notation and review basic graph theory and stability
definitions. In Section~\ref{sec:modeling} we formulate the rendezvous
control problem. The solution of the rendezvous control problem is
presented in Section~\ref{sec:soln_pcp_2stage}, together with an
intuitive description of its operation. The proof of the main theorem
is presented in Section~\ref{sec:proof}. Finally, in
Section~\ref{sec:conclusion} we make concluding remarks.  Lemmas and
claims related to the proof are in the appendix.

\section{Preliminaries}\label{sec:prelims}
\subsection{Notation}
We use interchangeably the notation $v=[v_1 \ \cdots \ v_n]^\top$ or
$(v_1,\ldots,v_n)$ for a column vector in $\Re^n$. We denote by $\one
\in \Re^m$ the vector $(1,\ldots,1)$.  If $v,w$ are vectors in
$\Re^2$, we denote by $v\cdot w:=v^\top w$ their Euclidean inner
product, and by $\|v\|:=(v \cdot v)^{1/2}$ the Euclidean norm of $v$.
%
%
If $\omega \in \Re$, we define
\[
\omega^\times  := \left[ \begin{array}{rr} 0   &  - \omega \\ 
\omega & 0
  \end{array}\right].
\]
Let $\{e_1,e_2\}$ denote the natural basis of $\Re^2$, $\SO(2):=\{M
\in \Re^{2 \times 2}: M^{-1} = M\trans, \det(M)=1\}$ and let ${\mathbb
  S}^1$ denote the unit circle. If $\Gamma$ is a closed subset of a
geodesically complete Riemannian manifold $\cX$, and $d: \cX \times
\cX \to [0,\infty)$ is a distance metric on $\cX$, we denote by
  $\|\chi\|_\Gamma:=\inf_{\psi \in \Gamma} d(\chi,\psi)$ the
  point-to-set distance of $\chi \in \cX$ to $\Gamma$. If
  $\varepsilon>0$, we let $B_\varepsilon(\Gamma) := \{\chi \in \cX :
  \|\chi\|_\Gamma < \varepsilon\}$ and by $\cN(\Gamma)$ we denote an
  open subset of $\cX$ containing $\Gamma$. If $A,B \subset \cX$ are
  two sets, denote by $A \backslash B$ the set-theoretic difference of
  $A$ and $B$. If $I = \{i_1,\ldots, i_n\}$ is an index set, the
  ordered list of elements $(x_{i_1},\ldots, x_{i_n})$ is denoted by
  $(x_j)_{j \in I}$.
  
Let $U, W$ be finite-dimensional vector spaces. A function $f: U
\rightarrow W$ is {\it homogeneous of degree $r$} if, for all
$\lambda>0$ and for all $x \in U$, $f(\lambda x)=\lambda^r f(x)$.  A
function $f: U \times V \rightarrow W$, $(x,y) \mapsto f(x,y)$, is
{\it homogeneous of degree $r$ with respect to $x$} if for all
$\lambda>0$ and for all $(x,y) \in U \times V$, $f(\lambda
x,y)=\lambda^r f(x,y)$.

\subsection{Graph Theory}\label{ssec:graph_theory}
We refer the reader to~\cite{godsil} for more details on the notions
reviewed in this section. We denote a {\em digraph} by
$\mathcal{G}=(\mathcal{V},\mathcal{E})$, where $\mathcal{V}$ is a set
of nodes labelled as $\{1,\dots,n\}$ and $\mathcal{E}$ is the set of
edges.   The set of {\em neighbors} of node $i$ is
$\cN_i:=\{ j \in {\cal V}: (i,j) \in \mathcal{E}\}$. 

Given positive numbers $a_{ij}>0$, $i,j \in \{1,\ldots,n\}$, the
associated {\em weighted Laplacian matrix} of $\cG$ is the matrix
$L:=D-A$, where $D$ is a diagonal matrix whose $i$-th diagonal entry
is the sum $\sum_{j \in \cN_i} a_{ij}$, and $A$ is the matrix whose
element $(A)_{ij}$ is $a_{ij}$ if $j \in \cN_i$, and $0$ otherwise.

 A {\it directed
  spanning tree} is a graph consisting of $n-1$ edges such that there
exists a unique directed path from a node, called the root, to every
other node. A {\it reverse directed spanning tree} is a graph which
becomes a directed spanning tree by reversing the directions of all
its edges. We identify the root of a reverse spanning tree with the
root of its associated spanning tree.  A graph ${\cal G}$ {\em contains a reverse
  directed spanning tree} if it has a subgraph which is a reverse
directed spanning tree.  

\begin{proposition}[\cite{beard,Lin}]\label{prop:connectivity}
The following conditions are equivalent for a digraph $\cG$:

\begin{enumerate}[(i)]

\item $\cG$ contains a reverse directed spanning tree.

\item For any set of positive gains $a_{ij}>0$, $i,j \in \{1,\ldots,
  n\}$ the associated weighted Laplacian matrix $L$ of $\cG$ has rank
  $n-1$, and $\Ker L = \Span\{\one\}$.
\end{enumerate}
\end{proposition}
A graph ${\cal G}=({\cal
  V},{\cal E})$ is \textit{strongly connected} if for any two nodes
$i,j \in \cal{V}$ there exists a path from $i$ to $j$.  A set of nodes
$S \subset {\cal V}$ is an \textit{isolated component} if it has no
outgoing edges, i.e., for any edge $(i,j) \in \cal{E}$, if $i \in S$
then $j \in S$. A graph ${\cal G'}=({\cal V'},{\cal E'})$ is a {\it
  subgraph} of ${\cal G}$ if ${\cal V'} \subset {\cal V}$ and ${\cal
  E'} \subset {\cal E}$. A subgraph ${\cal G'}$ is an {\em induced
  subgraph} of ${\cal G}$ if for any two vertices $i,j \in {\cal V'}$,
$(i,j) \in {\cal E'}$ if and only if $(i,j) \in {\cal E}$. A
\textit{strongly connected component} ${\cal G'}$ of ${\cal G}$ is a
maximal strongly connected induced subgraph of ${\cal G}$. In other
words, there does not exist any other strongly connected induced
subgraph of ${\cal G}$ containing ${\cal G'}$. Letting ${\cal
  G}_0=({\cal V}_0,{\cal E}_0), \dots, {\cal G}_r=({\cal V}_r,{\cal
  E}_r)$ be the strongly connected components of ${\cal G}$, the
\textit{condensation digraph} of ${\cal G}$, denoted ${\cal
  C(G)}=({\cal V_C({\cal G})},{\cal E_C({\cal G})})$, is defined as
follows. The \textit{vertex set} ${\cal V_C({\cal G})}$ is the set of
nodes $\{v_i\}_{i\in\{0,\dots,r\}}$ where the node $v_i$ is a {\em
  contraction} of the vertex set ${\cal V}_i$ of the $i$-th strongly
connected component ${\cal G}_i$. The \textit{edge set} ${\cal
  E_C({\cal G})}$ contains an edge $(v_i,v_j)$ if there exist vertices
$i' \in {\cal V}_i$ and $j' \in {\cal V}_j$ such that $(i',j') \in
{\cal E}$. The following properties of the condensation digraph are
found in~\cite{hatanaka2015}.
\begin{proposition}[\cite{hatanaka2015}]\label{prop:decomp}
Consider a graph ${\cal G}$ containing a reverse directed spanning
tree. The condensation ${\cal C(G)}$ satisfies the following
properties:
\begin{enumerate}[(i)]
\item ${\cal C(G)}$ is acyclic, i.e., there is no path in ${\cal
  C(G)}$ beginning and ending at the same node.

\item ${\cal C(G)}$ contains a reverse directed spanning tree ${\cal
  T}$ with a unique root $v_0 \in {\cal V_C({\cal G})}$.

\item There exists at least one vertex $v_i \in {\cal V_C({\cal G})}$
  such that $v_0$ is the only neighbor of $v_i$.
\end{enumerate}
\end{proposition} 
An example of a digraph ${\cal G}$ containing a reverse directed
spanning tree is shown in Figure~\ref{fig:graph}. The strongly
connected components are boxed. The resulting acyclic condensation
graph ${\cal C(G)}$ is shown in Figure~\ref{fig:graph2}. The
vertex $v_0$ in the figure is the unique root of the reverse directed
spanning tree in ${\cal C(G)}$.

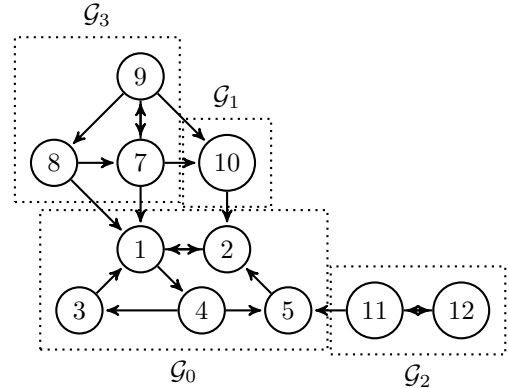
\begin{figure}[h]
\begin{center}
\begin{tikzpicture}[->,>=stealth',shorten >=1pt,auto,node distance=1.15cm,
  thick,main node/.style={circle,draw}]

  \node[main node] (1) {$1$};
  \node[main node] (2) [right of=1] {$2$};
\node[main node] (3) [below left of=1] {$3$};
\node[main node] (4) [below right of=1] {$4$};
\node[main node] (5) [right of=4] {$5$};

\node[main node] (7) [above of=1] {$7$};
\node[main node] (8) [left of=7] {$8$};
\node[main node] (9) [above of=7] {$9$};

\node[main node] (10) [right of=7] {$10$};

\node[main node] (11) [right of=5] {$11$};
\node[main node] (12) [right of=11] {$12$};
    
  \path 
  		(1) edge node[auto] {} (2)
		(2) edge node[auto] {} (1)
		(1) edge node[auto] {} (4)
		(4) edge node[auto] {} (3)
		(3) edge node[auto] {} (1)
		(4) edge node[auto] {} (5)
		(5) edge node[auto] {} (2)
		
		(7) edge node[auto] {} (1)
		(8) edge node[auto] {} (1)
		(8) edge node[auto] {} (7)
		(7) edge node[auto] {} (9)
		(9) edge node[auto] {} (7)
		(9) edge node[auto] {} (8)
		(7) edge node[auto] {} (10)
		(9) edge node[auto] {} (10)
		
		(10) edge node[auto] {} (2)
		
		(11) edge node[auto] {} (12)
		(12) edge node[auto] {} (11)
		(11) edge node[auto] {} (5);
		
\node[draw,inner sep=2mm,label=below:${\cal G}_0$,dotted,fit=(1) (5) (5) (3)] {};
\node[draw,inner sep=2mm,label=above:${\cal G}_3$,dotted,fit=(9) (7) (7) (8)] {};
\node[draw,inner sep=2mm,label=above:${\cal G}_1$,dotted,fit=(10) (10) (10) (10)] {};
\node[draw,inner sep=2mm,label=below:${\cal G}_2$,dotted,fit=(11) (11) (12) (11)] {};		
\end{tikzpicture}
\end{center}
\caption{Directed graph ${\cal G}$ containing a reverse directed spanning tree. The strongly connected components ${\cal G}_0, \dots, {\cal G}_3$ are boxed}
\label{fig:graph}
\end{figure}
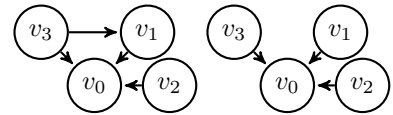
\begin{figure}[h]
\begin{center}
\begin{tikzpicture}[->,>=stealth',shorten >=1pt,auto,node distance=1cm,
  thick,main node/.style={circle,draw}]

  \node[main node] (1) {$v_0$};
  \node[main node] (2) [above right of=1] {$v_1$};
\node[main node] (3) [right of=1] {$v_2$};
\node[main node] (4) [above left of=1] {$v_3$};
    
  \path 
  		(2) edge node[auto] {} (1)
		(3) edge node[auto] {} (1)
		(4) edge node[auto] {} (1)
		(4) edge node[auto] {} (2);
		
\end{tikzpicture}
\begin{tikzpicture}[->,>=stealth',shorten >=1pt,auto,node distance=1cm,
  thick,main node/.style={circle,draw}]

  \node[main node] (1) {$v_0$};
  \node[main node] (2) [above right of=1] {$v_1$};
\node[main node] (3) [right of=1] {$v_2$};
\node[main node] (4) [above left of=1] {$v_3$};
    
  \path 
  		(2) edge node[auto] {} (1) [dashed]
  		(3) edge node[auto] {} (1)
		(4) edge node[auto] {} (1);
		
\end{tikzpicture}
\end{center}
\caption{Condensation digraph ${\cal C(G)}$ associated with the graph ${\cal G}$ in Figure \ref{fig:graph} (left) and reverse directed spanning tree contained in ${\cal C(G)}$ (right).}    
\label{fig:graph2}
\end{figure}

As in~\cite{hatanaka2015}, we define the vertex set ${\cal L}_k
\subset {\cal V}$ to be the union of those vertex sets ${\cal V}_i$
that correspond to vertices $v_i$ in the condensation digraph with the
property that the maximal path length from $v_i$ to the root $v_0$ is
equal to $k$.  By this definition, ${\cal L}_0:={\cal V}_0$. We let
${\cal L}_{-1}:= \varnothing$. Defining the vertex set $\bar{\cal
  L}_{k}:=\cup_{i=0}^{k}{\cal L}_i$,  by construction, the
neighbors of any vertex in ${\cal L}_k$ are contained in $\bar{\cal
  L}_{k-1}$. Therefore each node set
$\bar {\cal L}_k$ is isolated. For the example in Figure~\ref{fig:graph2}, we have
${\cal L}_0=\{ 1,2,3,4,5,6\}$, ${\cal L}_1=\{10\} \cup \{ 11, 12\}$
and ${\cal L}_2=\{7,8,9\}$.

\subsection{Stability Definitions}
The following stability definitions are taken from~\cite{ElHMag12}.
Let $\Sigma: \dot \chi = f(\chi)$ be a smooth dynamical system with
state space a geodesically complete Riemannian manifold $\cX$ with
Riemannian distance $d: \cX \times \cX \to [0,\infty)$, so that
  $(\cX,d)$ is a complete metric space. Let $\phi(t,\chi_0)$ denote
  the local phase flow of $\Sigma$. 
\begin{definition}
Consider a closed set $\Gamma
  \subset \cX$ that is positively invariant for $\Sigma$, i.e., for
  all $\chi_0 \in \Gamma$, $\phi(t,\chi_0) \in \Gamma$ for all $t>0$
  for which $\phi(t,\chi_0)$ is defined.

\begin{itemize}

\item $\Gamma$ is {\it stable} for $\Sigma$ if for any
$\varepsilon>0$, there exists a neighborhood $\cN(\Gamma) \subset \cX$
such that, for all $\chi_0 \in \cN(\Gamma)$, $\phi(t,\chi_0) \in
B_\varepsilon(\Gamma)$, for all $t >0$ for which $\phi(t,\chi_0)$ is
defined. 

\item $\Gamma$ is {\it attractive} for $\Sigma$ if there exists
  neighborhood $\cN(\Gamma) \subset \cX$ such that for all $\chi_0 \in
  \cN(\Gamma)$, $\lim_{t\to\infty} \|\phi(t,\chi_0)\|_\Gamma =0$.  The
     {\it domain of attraction of $\Gamma$} is the set $\{\chi_0 \in
     \cX: \lim_{t\to\infty} \| \phi(t,\chi_0)\|_\Gamma =0\}$.
     $\Gamma$ is {\it globally attractive} for $\Sigma$ if it is
     attractive with domain of attraction $\cX$.

\item $\Gamma$ is {\it locally asymptotically stable (LAS)} for
  $\Sigma$ if it is stable and attractive. The set $\Gamma$ is {\it
  globally asymptotically stable (GAS)} for $\Sigma$ if it is stable
  and globally attractive.  
\end{itemize}
\end{definition}
\begin{definition}
Let $\Gamma_1 \subset \Gamma_2$ be two subsets of $\cX$ that are
positively invariant for $\Sigma$. Assume that $\Gamma_1$ is compact
and $\Gamma_2$ is closed. 

\begin{itemize}

\item $\Gamma_1$ is {\it globally asymptotically stable relative to
  $\Gamma_2$} if it is GAS when initial conditions are restricted to
  lie in $\Gamma_2$.

\item $\Gamma_2$ is {\it locally stable near} $\Gamma_1$ if for all
  $c>0$ and all $\epsilon > 0$, there exists $\delta > 0$ such that
  for all $x_0 \in B_{\delta}(\Gamma_1)$ and all $t^\star>0$, if
  $\phi([0,t^\star],x_0) \subset B_c(\Gamma_1)$ then
  $\phi([0,t^\star],x_0) \subset B_{\epsilon}(\Gamma_2)$.

\item  $\Gamma_2$ is
  {\it locally attractive near} $\Gamma_1$ if there exists a
  neighbourhood $\cN(\Gamma_1)$ such that, for all $x_0 \in
  \cN(\Gamma_1), \ \| \phi(t,x_0)\|_{\Gamma_2} \rightarrow 0$ as $t
  \rightarrow \infty$.
\end{itemize}
\end{definition}
We present a reduction theorem used to derive our main result
 
\begin{theorem}[Reduction Theorem~\cite{ElHMag12,Seibert}]\label{thm:seibert-florio}
Let $\Gamma_1$ and $\Gamma_2$, $\Gamma_1 \subset \Gamma_2 \subset
\cX$, be two closed sets that are positively invariant for $\Sigma$,
and suppose $\Gamma_1$ is compact. Consider the following conditions:
(i) $\Gamma_1$ is LAS relative to $\Gamma_2$; (i') $\Gamma_1$ is GAS relative to $\Gamma_2$; (ii) $\Gamma_2$ is locally stable near $\Gamma_1$; (iii) $\Gamma_2$ is locally attractive near $\Gamma_1$; (iii)' $\Gamma_2$ is globally attractive; (iv) all trajectories of $\Sigma$ are bounded. 
%
%
%
%
%
%

Then, the following implications hold: 
(i) $\wedge$ (ii)
$\implies$ $\Gamma_1$ is stable; (i) $\wedge$ (ii)
$\wedge$ (iii) $\iff$ $\Gamma_1$ is LAS; (i)' $\wedge$ (ii) $\wedge$ (iii)' $\wedge$ (iv) $\iff$
  $\Gamma_1$ is GAS.
%
%
%
%
%
%
\end{theorem}

\section{Rendezvous Control Problem}\label{sec:modeling}

Consider a group of $n$ kinematic unicycles.  Let $\cI=\{i_x,i_y\}$ be
an inertial frame in three-dimensional space and consider the $i$-the
unicycle in Figure~\ref{fig:vehicle_detailed}.  Fix a body frame
$\cB_{i}=\{b_{ix},b_{iy}\}$ to the unicycle, where $b_{ix}$ is the
heading axis, and denote by $x_i \in \Re^2$ the position of the
unicycle in the coordinates of frame $\cI$. The unicycle's attitude is
represented by a rotation matrix $R_i$ whose columns are the
coordinate representations of $b_{ix}$ and $b_{iy}$ in frame $\cI$.
Letting $\theta_i \in \mathbb{S}^1$ be the angle between vectors $i_x$
and $b_{ix}$, we have
\[
R_i=\begin{bmatrix}
\cos \theta_i& -\sin \theta_i\\
\sin \theta_i& \cos \theta_i\\
\end{bmatrix}.
\]
The angular speed of robot $i$ is denoted by $\omega_i$. The unicycle
dynamics are given by,
\begin{align}
& \dot{x}_{i} = u_iR_ie_1 \label{eq:vehicle:translational}\\ 
& \dot R_{i} = R_i(\omega_i)^\times, \quad i=1,\dots,
  n. \label{eq:vehicle:rotational:so3}
\end{align}
In what follows, we refer to
system~\eqref{eq:vehicle:translational}-\eqref{eq:vehicle:rotational:so3}
as $\Sigma_i$.  Its control inputs are the linear speed $u_i$ and
angular speed $\omega_i$. The relative displacement of robot $j$ with
respect to robot $i$ is $x_{ij}:=x_j - x_i$. If $v \in \Re^2$ is the
coordinate representation of a vector in frame $\cI$, then we denote
by $v^i:=R_i^{-1} v$ the coordinate representation of $v$ in body
frame $\cB_i$.

We define the {\em sensor digraph}
$\mathcal{G}=(\mathcal{V},\mathcal{E})$, where each node represents a
robot, and an edge from node $i$ to node $j$ indicates that robot $i$
can sense robot $j$. We assume that ${\cal G}$ has no self-loops and
is time-invariant. Given a node $i$, its set of neighbors $\cN_i$
represents the set of vehicles that robot $i$ can sense. If $j \in
{\cal N}_i$, then we say that robot $j$ is a {\em neighbour} of robot
$i$. If this is the case, then robot $i$ can sense the relative
displacement of robot $j$ in its own body frame, i.e., the quantity
$x_{ij}^i$.  Define the vector $y_i:=(x_{ij})_{j \in {\cal N}_i}$.
The relative displacements available to robot $i$ are contained in the
vector $y_i^i:=(x^{i}_{ij})_{j \in {\cal N}_i}$. A {\it local and
  distributed feedback} $(u_i,\omega_i)$ for robot $i$ is a locally
Lipschitz function of $y_i^i$. We define the {\em rendezvous manifold}
\begin{equation} \label{eq:rendez_set}
\Gamma := \left\{(x_i,R_i)_{i \in \{1, \dots, n\}} \in \Re^{2n} \times
\SO(2)^n: x_{ij}=0, \ \forall\, i,j \right\}.
\end{equation}

We are now ready to state the rendezvous control problem.

\vspace*{2ex}
\noindent {\em Rendezvous Control Problem:} For
system~\eqref{eq:vehicle:translational}-\eqref{eq:vehicle:rotational:so3}
with sensor digraph ${\cal G}$, find local and distributed feedbacks
$(u_i,\omega_i)_{i \in \{1,\ldots,n\}}$ that globally asymptotically
stabilize the rendezvous manifold $\Gamma$. \hfill $\triangle$

\section{Solution of the Rendezvous Control Problem}\label{sec:soln_pcp_2stage}

\begin{figure}[h]
\vspace*{0.1in}
\centerline{\includegraphics[width=0.45\textwidth]{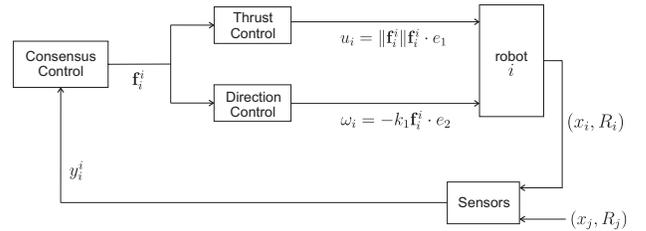}}
\caption{Block diagram of the rendezvous control system for robot
  $i$.}
\label{fig:block_diagram}
\end{figure}
In this section we present the solution of the rendezvous control
problem. Consider the function,
\begin{equation}\label{eqn:linear}
\wdes_i(y_i):=\sum_{j\in {\cal N}_i}a_{ij}x_{ij},
\end{equation}
$i=1, \dots, n$ with $a_{ij} > 0$. The function $\wdes_i(y_i)$ is a
standard linear consensus controller for single integrator systems
\cite{beard,moreau,olfati}. We use $\wdes(y_i)$ to construct the
feedbacks

\begin{equation} \label{eq:CCP}
\begin{aligned}
&u_i =\|\wdes_i(y_i^i)\|\wdes_i(y_i^i) \cdot e_1, \\
&\omega_i=-\kq \wdes_i(y_i^i) \cdot e_2, \ i=1, \dots, n.
\end{aligned}
\end{equation}

The result below states that for sufficiently large $\kq$, the
feedbacks in~\eqref{eq:CCP} solve the rendezvous control problem if the
network of unicycles has a sensor digraph containing a reverse
directed spanning tree. 

\begin{theorem}\label{thm:main3}
The rendezvous control problem is solvable for
system~\eqref{eq:vehicle:translational}-\eqref{eq:vehicle:rotational:so3}
if, and only if, the sensor graph ${\cal G}$ contains a reverse
directed spanning tree, in which case a solution is as follows.  There
exists $\kq^\star > 0$ such that for any $\kq>\kq^\star$
feedback~\eqref{eq:CCP} with $\wdes_i(y_i)$ in~\eqref{eqn:linear}
solves the rendezvous control problem.
\end{theorem}

The necessity portion of Theorem~\ref{thm:main3} was proved
in~\cite{Lin}. The sufficiency part, namely the fact that the
feedback~\eqref{eq:CCP} solves the rendezvous control problem, is
proved in Section~\ref{sec:proof}. 

\begin{figure}[h]
\centerline{\includegraphics[width=0.20\textwidth]{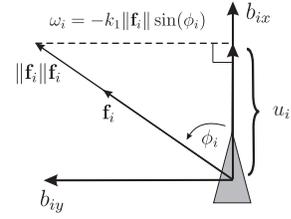}}
\caption{Illustration of the control inputs $u_i$ and $\omega_i$ in~\eqref{eq:CCP}.}
\label{fig:vehicle_detailed}
\end{figure}

The proposed control architecture is illustrated in the block diagram
of Figure~\ref{fig:block_diagram}. There are two nested loops. The
outer loop treats each robot as a single-integrator driven by the
linear consensus controller,
\begin{equation} \label{eq:dcg}
\dot x_i=\wdes_i(y_i), \ i=1,\dots,n.
\end{equation}
The set $\left\{(x_i)_{i \in \{1, \dots, n\}} \in \Re^{2n}: x_{ij}=0,
\ \forall i,j \right\}$ is globally asymptotically stable
for~\eqref{eq:dcg} if the sensing graph has a reverse directed
spanning tree \cite{moreau}. The signal $\wdes_i(y_i(t))$ is computed
in the body frame $\cB_i$, and used as a reference signal for the
inner-loop thrust and direction controllers that assign the unicycle
control inputs in~\eqref{eq:CCP}. The intuition behind these
controllers is shown in Figure~\ref{fig:vehicle_detailed}. The
speed input $u_i$ is the dot product
$u_i=\|\wdes_i(y_i^i)\|\wdes_i(y_i^i) \cdot e_1$. This is the
projection of the reference $\|\wdes_i(y_i)\|\wdes_i(y_i)$ onto the
heading axis $b_{ix}$ of robot $i$. The angular speed, on the other
hand, is proportional to the dot product between the reference
$\wdes_i(y_i)$ and the second body axis $b_{iy}$.  In
Figure~\ref{fig:vehicle_detailed}, one can see that
$\omega_i=-k_1\|\wdes_i\|\sin(\phi_i)$ acts to reduce the angle
$\phi_i$ between $b_{ix}$ and $\wdes_i(y_i)$ with a rate proportional
to the magnitude of $\wdes_i$. Together, these control inputs drive
the robot velocity $u_i b_{ix}$ approximately to the reference
$\|\wdes_i(y_i)\|\wdes_i(y_i)$. The convergence is approximate because
the control inputs do not depend on the time derivative of
$\wdes_i$. It is the difference in angle between $u_i b_{ix}$ and
$\|\wdes_i(y_i)\|\wdes_i(y_i)$ as opposed to the difference in
magnitude that is important for obtaining rendezvous. Since
$\|\wdes_i(y_i)\|\wdes_i(y_i)$ is homogeneous of degree two, as the
robots approach consensus, $\omega_i$ converges to zero slower than
$u_i$. This allows $\omega_i$ to exert sufficient control authority
even as the robots converge to consensus, closing the gap between the
vectors $u_i b_{ix}$ and $\|\wdes_i(y_i)\|\wdes_i(y_i)$.

\subsection{Simulation Results}
\begin{figure*}[t]
\centerline{\includegraphics[width=0.9\textwidth]{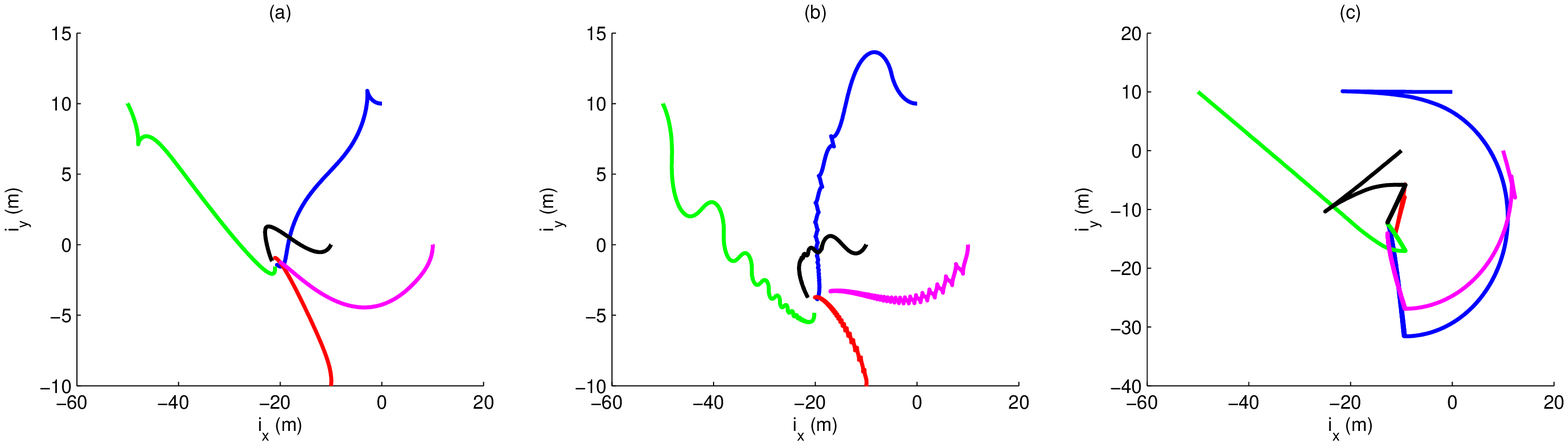}}
\caption{Rendezvous control simulation for: (a) proposed feedback in~\eqref{eq:CCP}, (b) feedback in~\cite{Lin}, and (c) feedback in~\cite{dimarogonas}}
\label{fig:sim1}
\end{figure*}

We consider a group of five robots with sensor digraph in Figure~\ref{fig:simgraph}. For the feedback in~\eqref{eq:CCP}, we pick $a_{ij}=0.05$ for all $j \in {\cal N}_i$. The control gain $\kq$ is chosen to be $\kq=1$. The initial conditions of the robots are shown in Table~\ref{tab:sim}. The simulation is presented in Figure~\ref{fig:sim1}(a). The proposed feedback has practical advantages over the time-varying feedback in~\cite{Lin} and the discontinuous feedback in~\cite{dimarogonas} whose simulation results are shown in Figure~\ref{fig:sim1}(b) and Figure~\ref{fig:sim1}(c) respectively with the same initial conditions in Table~\ref{tab:sim} and sensing graph in Figure~\ref{fig:simgraph}. The proposed feedback induces a more natural behaviour in the ensemble of unicycles. The feedback in~\cite{Lin} makes the unicycle ``wiggle'' indefinitely, a behaviour which would be unacceptable in practice. The feedback in~\cite{dimarogonas} induces instantaneous changes in direction that are impossible to achieve with realistic implementations.
\begin{figure}[h]
\begin{minipage}{0.15\textwidth}
\begin{center}
\begin{tikzpicture}[->,>=stealth',shorten >=0.5pt,auto,node distance=1cm,
  thick,main node/.style={circle,draw}]

  \node[main node] (3) {$3$};
  \node[main node] (1) [right of=3] {$1$};
  \node[main node] (2) [above of=3] {$2$};
  \node[main node] (4) [left of=3] {$4$};
  \node[main node] (5) [above of=1] {$5$};
    
  \path 
  		(1) edge node[auto] {} (2)
        (1) edge node[auto] {} (3)
        (2) edge node[auto] {} (4)
        (3) edge node[auto] {} (2)
        (4) edge node[auto] {} (3)
        (5) edge node[auto] {} (2);
\end{tikzpicture}
\end{center}
\caption{Sensor digraph used in the simulation results.}
\label{fig:simgraph}
\end{minipage}
\begin{minipage}{0.35\textwidth}
\captionof{table}{Simulation Initial Conditions}
\centering
\begin{tabular}{c c c c}
\hline\hline
Vehicle $i$ & $x_i(0)$ (m) & $\theta_i(0)$ (rad)\\
\hline
1 & $(0,10)$ & $0$\\
2 & $(-10,-10)$ & $2\pi/5$\\
3 & $(-50,10)$ & $4\pi/5$ \\
4 & $(-10,0)$ & $6\pi/5$ \\
5 & $(10,0)$ & $8\pi/5$ \\[1ex]
\hline
\end{tabular}
\label{tab:sim}
\end{minipage}
\end{figure}
%
%
%

\section{Proof of Theorem~\ref{thm:main3}}\label{sec:proof}

This section presents the sufficiency proof of
Theorem~\ref{thm:main3}. The necessity was proved in~\cite{Lin}.  The
key tool in our proof is the condensation graph and the isolated node
sets $\bar {\cal L}_k$ defined in Section~\ref{ssec:graph_theory}.
The same tool was employed in~\cite{hatanaka2015} for pose
synchronization (synchronization of positions and attitudes) of fully
actuated vehicles.

The dynamics of unicycles associated with an isolated node set $\bar
\cL_k$ are independent of the nodes outside of this set because, for
any robot $i \in \bar \cL_k$, the feedbacks $u_i$ and $\omega_i$
in~\eqref{eqn:linear}, \eqref{eq:CCP} depend only on states of robots
within $\bar \cL_k$.  Therefore, the dynamics of the collection of
unicycles in $\bar \cL_k$,
\begin{align}
\begin{split} 
	\dot{x}_{i} = u_iR_ie_1 \label{eq:rel_coords_trans_iso}\\
\end{split}
\\
\begin{split} 
	\dot R_{i} = R_i(\omega_i)^\times, \quad i \in \bar \cL_k \label{eq:rel_coords_rot_iso}\\
\end{split}
\end{align}
define an autonomous dynamical system.  Henceforth, the dynamics
in~\eqref{eq:rel_coords_trans_iso},~\eqref{eq:rel_coords_rot_iso} are
denoted by $\Sigma_{\bar \cL_k}$ and we
define the {\em reduced rendezvous manifold} $\Gamma_{\bar
  \cL_k}:=\left\{(x_i,R_i)_{i \in \bar \cL_k}: x_{ij}=0, \ \forall i,j
\in \bar \cL_k \right\}$.

Recall from Section~\ref{ssec:graph_theory} that the set $\bar
\cL_{-1}$ is empty, which implies that the set $\Gamma_{\bar
  \cL_{-1}}$ is also empty. We adopt the convention that $\Gamma_{\bar
  \cL_{-1}}$ is GAS for $\Sigma_{\bar \cL_{-1}}$.

The proof of Theorem~\ref{thm:main3} relies on an induction argument
on the node sets $\bar \cL_k$. Key in the induction argument is the
next result stating that if the vehicles in $\bar {\cal L}_{k-1}$
achieve rendezvous, then so do the vehicles in $\bar {\cal L}_k$.

\begin{proposition}\label{thm:main2}
Consider system~\eqref{eq:vehicle:translational},
~\eqref{eq:vehicle:rotational:so3} and assume that the sensor graph
${\cal G}$ contains a reverse directed spanning tree. Let $u_i$ and
$\omega_i$ be as in~\eqref{eq:CCP} with $\wdes_i(y_i)$ as
in~\eqref{eqn:linear}. Suppose that, for some integer $k \geq 0$, the
set $\Gamma_{\bar {\cal L}_{k-1}}$ is globally asymptotically stable
for the dynamics $\Sigma_{\bar {\cal L}_{k-1}}$. There exists
$\kq^\star > 0$ such that choosing $\kq>\kq^\star$ in~\eqref{eq:CCP},
implies $\Gamma_{\bar {\cal L}_k}$ is globally asymptotically stable
for the dynamics $\Sigma_{\bar {\cal L}_k}$.
\end{proposition}

In Section~\ref{ssec:main_proof}, we use the above proposition to
prove Theorem~\ref{thm:main3}, and in Section~\ref{ssec:prelim_proof}
we prove Proposition~\ref{thm:main2}.

In the special case when ${\cal G}$ is strongly connected, we have
$\bar {\cal L}_0 = {\cal V}$. Since, by definition, $\bar {\cal
  L}_{-1} = \varnothing$,  the set $\Gamma_{\bar {\cal
    L}_{-1}}$ is GAS for $\Sigma_{\bar {\cal L}_{-1}}$, and
Proposition~\ref{thm:main2} yields the following corollary.

\begin{corollary}\label{thm:main}
Consider system~\eqref{eq:vehicle:translational},
~\eqref{eq:vehicle:rotational:so3} and assume that the sensor graph
${\cal G}$ is strongly connected. Let $u_i$ and $\omega_i$ be as
in~\eqref{eq:CCP} with $\wdes_i(y_i)$ as in~\eqref{eqn:linear}. There
exists $\kq^\star > 0$ such that choosing $\kq>\kq^\star$ solves the
rendezvous control problem.
\end{corollary}

\subsection{Proof of Theorem~\ref{thm:main3}}\label{ssec:main_proof}
To begin with, the feedback in~\eqref{eq:CCP} is local and distributed because it is
a smooth function of $y_i^i$ only. Consider a graph $\cal{G}=(\cal{V},\cal{E})$ containing a reverse directed spanning tree and the node sets ${\cal L}_k$ and $\bar {\cal L}_{k}$ defined in Section~\ref{ssec:graph_theory}. By construction, the node sets $\bar {\cal L}_{k}$ are isolated, the subgraph $({\cal V}_0,{\cal E}_0)$ is strongly connected, and $\bar {\cal L}_0 = {\cal L}_0 = {\cal V}_0$. 

The proof is by induction. Since the subgraph $(\bar {\cal L}_0,{\cal E}_0)$ is strongly connected, by Corollary~\ref{thm:main}, there exists $l_0$ such that choosing $k_1>l_0$ makes the set $\Gamma_{\bar {\cal L}_0}$ globally asymptotically stable for system $\Sigma_{\bar {\cal L}_0}$. 

Now consider $\bar {\cal L}_k$ and suppose the reduced rendezvous
manifold $\Gamma_{\bar {\cal L}_{k-1}}$ is globally asymptotically
stable for system $\Sigma_{\bar {\cal L}_{k-1}}$. It holds from
Proposition~\ref{thm:main2} that there exists $l_k$ such that choosing
$k_1>l_k$ makes the isolated node set $\Gamma_{\bar {\cal L}_k}$
globally asymptotically stable for system $\Sigma_{\bar {\cal L}_k}$. By part (ii) of Proposition~\ref{prop:decomp}, ${\cal C(G)}$ contains a reverse directed spanning tree, so there is a path from every node of ${\cal C(G)}$ to the unique root of ${\cal C(G)}$. By part (i) of the same proposition, ${\cal C(G)}$ is acyclic, which implies that the paths connecting the nodes of ${\cal C(G)}$ to the unique root of ${\cal C(G)}$ have a maximum length, $k^\star$. Recall that, by definition, $\bar \cL_{k^\star} = \sum_{i=1}^{k^\star} \cL_i$ is the union of those strongly connected components ${\cal V}_i$ of ${\cal V}$ that are associated with nodes $v_i$ of the condensation digraph ${\cal C(G)}$ with the property that the maximum path length from $v_i$ to the root $v_0$ is $\leq k^\star$. As we argued earlier, the set of such nodes $v_i$ equals the entire condensation digraph, implying that  $\bar \cL_{k^\star} = {\cal V}$. Let
$k_1^\star>\max\{l_0,\dots,l_{k^\star}\}$. By induction, it must hold
that choosing $k_1>k_1^\star$ makes $\Gamma_{\bar {\cal
    L}_{k^\star}}=\Gamma$ globally asymptotically stable for system
$\Sigma_{\bar {\cal L}_{k^\star}}=\Sigma_{\cal{V}}=\Sigma$. We
conclude that $\Gamma$ is globally asymptotically stable. \qed

\subsection{Proof of Proposition~\ref{thm:main2}}\label{ssec:prelim_proof}
We denote $A:=\bar {\cal L}_{k-1}$ and $B:={\cal L}_k$ and therefore $\bar {\cal L}_k= A \cup B$. By assumption, $\Gamma_A$ is globally asymptotically stable for the dynamics $\Sigma_A$ and the graph associated to the nodes in $B$ is strongly connected. We need to show that $\Gamma_{A \cup B}$ is globally asymptotically stable for the dynamics $\Sigma_{A \cup B}$. The proof relies on the following coordinate transformation.

\subsubsection{Coordinate Transformation}
For notational convenience, we collect the position vectors $x_i$ and
rotation matrices $R_i$ into variables $x:=(x_1,\ldots,x_n)$ and
$R:=(R_1,\ldots, R_n)$. We define the spaces
$
\sX :=\Re^{2n}, \ \sR :=\SO(2) \times \cdots \times \SO(2) \ (n \text{
  times}),
$
so that $x \in \sX$ and $R \in \sR$.  For each $i
\in \{1,\ldots,n\}$, define 
\begin{equation}\label{eq:Chi}
\Chi_i := \wdes_i(y_i) / A_i,  
\end{equation}
where $A_i :=\sum_{j \in \cN_i} a_{ij}$, and let $\Chi:=(\Chi_1,\ldots,
\Chi_n)$. We may express $\Chi$ as
\[
\Chi = \diag(1/A_1,\cdots,1/A_n) (L \otimes I_2) x.
\]
In the above, $\diag(\ldots)$ is the diagonal matrix with diagonal
elements inside the parenthesis; $L$ is the weighted Laplacian matrix
of the sensor digraph associated with the gains $a_{ij}$; finally,
$\otimes$ denotes the Kronecker product of matrices.  Since the sensor
digraph contains a reverse directed spanning tree, by
Proposition~\ref{prop:connectivity} the matrix $L\otimes I_2$ has rank
$2(n-1)$, and $\Ker (L\otimes I_2) = \Span\{\one \otimes e_1,\one \otimes e_2\}$ with $\one \in \Re^{n}$. Let $\bar x :=[I_2
  \ \cdots \ I_2]x = \sum_i x_i$, then the linear map $T: \sX \to \sX
\times \Re^2$, $x \mapsto (\Chi,\bar x)$ is an isomorphism onto its
image.  Under the action of $T$, the subspace $\{x \in \sX: x_1 =
\cdots =x_n\}$ is mapped isomorphically onto the subspace
$\{(\Chi,\bar x) \in \image T: \Chi=0\}$. Since the feedbacks
in~\eqref{eqn:linear}-\eqref{eq:CCP} are local and distributed, it can
be seen that the dynamics of the closed-loop unicycles in $(\Chi,\bar
x, R)$ coordinates are independent of $\bar x$. Moreover, as we have
seen, in these coordinates the control specification is the global
stabilization of $\{(\Chi,\bar x,R) \in \sX \times \Re^2 \times \sR:
\Chi=0\}$, a set whose description is independent of $\bar x$. In
light of these considerations, for the stability analysis we may drop
the variable $\bar x$, and show that the set
$
\hat \Gamma:=\{(\Chi,R) \in \sX \times \sR: \Chi =0\}
$
is GAS for the $(\Chi,R)$ dynamics. 

From here on we will use the hat notation to refer to quantities
represented in $(\Chi,R)$ coordinates. Denote $\gdes_i(y_i)
:=\|\wdes_i(y_i)\| \wdes_i(y_i)$. Using~\eqref{eq:Chi}, the functions
$\wdes_i$ and $\gdes_i$ and their body frame representations are given
in $(\Chi,R)$ coordinates by
\begin{equation}\label{eq:g_and_OO}
\begin{aligned}
&\hat \wdes_i(\Chi_i)=\p \Chi_i, \
\hat \gdes_i(\Chi_i)=\p^2 \|\Chi_i\|\Chi_i \\
&\hat \wdes_i^i(\Chi_i,R_i)=\p R_i^{-1} \Chi_i, \ 
\hat \gdes_i^i(\Chi_i,R_i)=\p^2 R_i^{-1}
\|\Chi_i\|\Chi_i,
\end{aligned}
\end{equation}
and we can use these expressions to rewrite the
feedback~\eqref{eq:CCP} in new coordinates as
$
u_i = \hat \gdes_i(\Chi_i,R_i) \cdot e_1, \ \omega_i = -k_1 \hat \wdes_i(\Chi_i,R_i) \cdot e_2.
$
We remark that $\hat \wdes_i$ and $\hat \wdes_i^i$ are homogeneous of
degree one with respect to $\Chi_i$. Similarly, $\hat \gdes_i$ and
$\hat \gdes_i^i$ are homogeneous of degree two with respect to
$\Chi_i$.  The closed-loop unicycle dynamics in $(\Chi,R)$ coordinates
are given by
\begin{align}
\begin{split}
	\dot \Chi_i &= \frac{\sum_{j \in N_i}a_{ij}((\hat \gdes_j^j
          \cdot e_1)R_je_1-(\hat \gdes_i^i \cdot e_1)R_ie_1)}{\p}, \label{eq:rel_coords_trans}\\
\end{split}
\\
\begin{split} 
	\dot R_{i} = R_i(-\kq \hat \wdes_i^i \cdot
        e_2)^\times. \label{eq:rel_coords_rot}\\
\end{split}
\end{align}
We will refer to
system~\eqref{eq:rel_coords_trans}-\eqref{eq:rel_coords_rot} as $\hat
\Sigma_i$.

In analogy with what we did earlier, for a set of nodes $S \subset
{\cal V}$ we let $\Chi_S:=(\Chi_i)_{i\in S} \in \sX_S$ and
$R_S:=(R_i)_{i\in S} \in \sR_S$. Moreover, if $S$ is an isolated node
set, the systems $\hat \Sigma_i, i \in S$ determine an autonomous
dynamical system which we denote by $\hat \Sigma_S$. We also denote
the reduced rendezvous manifold by
$
\hat\Gamma_S:=\left\{(\Chi_S,R_S) \in \sX_S \times \sR_S: \Chi_S=0\right\}.
$
In new coordinates, it needs to be shown that the set $\hat \Gamma_{A \cup B}$ is globally asymptotically stable for the dynamics $\hat \Sigma_{A \cup B}$ under the assumption that $\hat \Gamma_A$ is globally asymptotically stable for the dynamics $\hat \Sigma_A$.
\subsubsection{Stability analysis}

Let 
\begin{equation}\label{eq:lyaps}
\begin{aligned}
V(\Chi_B)&=\sum_{i \in B} \gamma_i \Chi_i^\top \Chi_i \\
W_{\tran}(\Chi_B) &= \sqrt{V(\Chi_B)} \\
W_{\rot}(\Chi_B,R_B) &=\sum_{i \in B} \hat \wdes_i^i (\Chi_i,R_i) \cdot e_1,
\end{aligned}
\end{equation}
where $\gamma_i>0$ are gains that will be defined later. 
 Consider
the function $W:\sX_B \times \sR_B \rightarrow \Re$ defined as
\begin{equation}\label{eq:W}
 W(\Chi_B,R_B)=\alpha W_{\tran}(\Chi_B)+ W_{\rot}(\Chi_B,R_B),
\end{equation}
where  $\alpha>0$ is a design parameter.

The next two lemmas are used in the subsequent analysis.

\begin{lemma}\label{lem:saturation}
Consider the continuous function $W(\Chi_B,R_B)$ defined in~\eqref{eq:W}. There
exists $\alpha^\star>0$ such that, for all $\alpha > 2\alpha^\star$,
the following properties hold:
\begin{enumerate}[(i)]
\item $W \geq 0$ and $W^{-1}(0) = \{(\Chi_B,R_B): \Chi_B=0\}$.
\item For all $c>0$, the sublevel set $W_c :=\{(\Chi_B,R_B):
  W(\Chi_B,R_B) \leq c\}$ is compact.
\item $\alpha^\star \sqrt{V(\Chi_B)} <W(\Chi_B,R_B)< 2\alpha \sqrt{V(\Chi_B)}$.
\end{enumerate}
\end{lemma}
The proof is in the appendix. From now on assume $\alpha > 2\alpha^\star$. 

\begin{lemma}\label{lem:inequalities}
Consider system~\eqref{eq:rel_coords_trans},~\eqref{eq:rel_coords_rot}.
There exist gains $\gamma_i$ in~\eqref{eq:lyaps} and $\kq^\star>0$ such that choosing $\kq>\kq^\star$ implies 
\begin{equation}\label{eq:ineq}
\frac{d}{dt} W(\Chi_B,R_B) \leq -\sigma V(\Chi_B) + \Phi(\Chi_A,R), \ \sigma>0,
\end{equation}
where $\Phi(\Chi_A,R)$ is continuous with respect to its arguments and $\Phi(0,R)=0$.
\end{lemma}
The proof of Lemma~\ref{lem:inequalities} is presented in the appendix.

We will now show that choosing $\kq>\kq^\star$ implies $\hat \Gamma_{A
  \cup B}$ is globally asymptotically stable for $\hat
\Sigma_{A \cup B}$. The proof will make use of the reduction theorem
(Theorem~\ref{thm:seibert-florio}).  We will first show that all
solutions of the closed-loop system are bounded. The rotation matrices
live in a compact set, therefore we only need to show that the states
$\Chi_{A \cup B}=(\Chi_i)_{i \in A \cup B}$ are bounded. Since $A$ is
isolated, $\hat \Sigma_A$ is an autonomous subsystem and by
assumption, $\hat \Gamma_A=\{(\Chi_A,R_A) \in \sX_A \times \sR_A:
\Chi_A=0\}$ (compact), is globally asymptotically stable. Therefore,
$\Chi_A$ is bounded. From the inequality $W(\Chi_B,R_B)\geq
\alpha^\star \sqrt{V(\Chi_B)}$ in part (iii) of
Lemma~\ref{lem:saturation}, to show boundedness of $V(\Chi_B)$, it
suffices to show that $W(\Chi_B,R_B)$ is bounded. Boundedness of
$V(\Chi_B)$, in turn, implies boundedness of $\Chi_B$. From the bound
on the derivative of $W$ in~\eqref{eq:ineq}, and by
Lemma~\ref{lem:saturation} we obtain
\[
\frac{d}{dt} W(\Chi_B,R_B) \leq -\frac{\sigma W(\Chi_B,R_B)^2}{ (2\alpha)^2} + \Phi(\Chi_A,R), \ \sigma>0.
\]
Since $\Chi_A$ is bounded and $R \in \sR$ lies on a compact set, it holds that $\Phi(\Chi_A,R)$ is bounded and therefore $W$ is bounded, which implies that $\Chi_B$ is bounded. Therefore $\Chi_{A \cup B}$ is bounded, as claimed. Now define the set,
$
\hat \Lambda:=\{(\Chi_{A \cup B},R_{A \cup B}) \in \sX_{A \cup B} \times \sR_{A \cup B}: \Chi_A=0\}.
$
Since the set $\hat \Gamma_A$ is globally asymptotically stable for system $\hat \Sigma_A$ and $\Chi_{A \cup B}$ is bounded, it holds that $\hat \Lambda$ is globally asymptotically stable for $\hat \Sigma_{A \cup B}$.  

To show that the set $\hat \Gamma_{A \cup B}$, which is compact, is globally asymptotically stable for the system $\hat \Sigma_{A \cup B}$, it suffices to show that $\hat \Gamma_{A \cup B}$ is globally asymptotically stable relative to $\hat \Lambda$. On the set $\hat \Lambda$, $\Phi(\Chi_A,R)$ is equal to zero and the derivative of $W$ is therefore given by
$
\frac{d}{dt} W(\Chi_B,R_B) \leq -\frac{\sigma W(\Chi_B,R_B)^2}{(2 \alpha)^2}, \ \sigma>0.
$
By Lemma~\ref{lem:saturation}, all level sets of $W(\Chi_B,R_B)$ are compact and $W^{-1}(0) = \{(\Chi_B,R_B): \Chi_B=0\}$. This implies $\hat \Gamma_{A \cup B}$ is globally asymptotically stable relative to the set $\hat \Lambda$. By Theorem~\ref{thm:seibert-florio}, $\hat \Gamma_{A \cup B}$ is globally asymptotically stable for $\hat \Sigma_{A \cup B}$. This completes the proof.

\section{Conclusion}\label{sec:conclusion}
We have presented the first solution to the rendezvous control problem for a group of kinematic unicycles on the plane using continuous, time-independent feedback that is local and distributed. The solution assumes a fixed sensing digraph that contains a reverse-directed spanning tree. The control methodology is based on a control structure made of two nested loops. An outer loop produces a standard feedback for concensus of single integrators which becomes reference to an inner loop assigning the unicycle control inputs that rely only on onboard measurements. Information of the unicycle's relative orientations is not required.

\bibliographystyle{IEEEtran} \bibliography{sub}
\appendix

Throughout this appendix we will make use of functions $\mu_i$ and $\mu$ defined as follows. 
Recall that  $V(\Chi_B)$ is positive definite.  Define the functions $\mu: \sX_B \backslash 0 \rightarrow \mu(\sX_B\backslash 0)$, $\mu(\Chi_B):=\Chi_B/\sqrt{V(\Chi_B)}$, and $\mu_i: \sX_B\backslash 0 \rightarrow \mu_i(\sX_B\backslash 0)$, $\mu_i(\Chi_B):=\Chi_i/\sqrt{V(\Chi_B)}, \ i \in B$. Since the numerator and denominator are both homogeneous of degree one, these functions are both homogeneous of degree zero with respect to $\Chi_B$. Therefore, the images satisfy $\mu(\sX_B\backslash 0)=\mu({\mathbb S}^k)$ and $\mu_i(\sX_B\backslash 0)=\mu_i({\mathbb S}^k)$, where ${\mathbb S}^k$ is the unit sphere in $\sX_B$. Since $\mu$ and $\mu_i$ are continuous functions and ${\mathbb S}^k$ is a compact set, the images $\mu(\sX_B\backslash 0)$ and $\mu_i(\sX_B\backslash 0)$ are compact sets.

\subsection{Proof of Lemma~\ref{lem:saturation}}

Recall the definition of $W(\Chi_B,R_B)$, 
\[
\begin{aligned}
& W = \alpha \sqrt{V(\Chi_B)} +
\sum_{i \in B} \hat\wdes_i^i(\Chi_i,R_i)\cdot e_1 \\
& \quad =\sqrt{V(\Chi_B)} \left( \alpha + \frac{ \sum_{i \in B} \hat\wdes_i^i(\Chi_i,R_i)\cdot
  e_1 }{\sqrt{V(\Chi_B)}} \right).
\end{aligned}
\]
Using the fact that $\hat\wdes_i^i(\Chi_i,R_i)$ is homogeneous with respect to
its first argument, we have
$
W=\sqrt{V(\Chi_B)} \left( \alpha + \sum_{i \in B}
\hat\wdes_i^i\left(\mu_i(\Chi_B),R_i\right)\cdot e_1 \right).
$
Since $\hat\wdes_i^i$ is continuous, $\mu_i(\Chi_B)$ is bounded, and $R_B \in
\sR_B$, a compact set, it follows that the function $\sum_{i\in B} \left|
\hat\wdes_i^i\left(\mu_i(\Chi_B),R_i\right)\cdot e_3 \right|$ has a bounded
supremum. Accordingly, let 
$
\alpha^\star = \sup_{(\Chi_B,R_B)\in \sX_B \times \sR_B} \sum_{i \in B} \left|
\hat\wdes_i^i\left(\mu_i(\Chi_B),R_i\right)\cdot e_1 \right|.
$
For all $\alpha > 2\alpha^\star$, we have
$
W(\Chi_B,R_B) \geq \underline{W}(\Chi_B,R_B):= \alpha^\star \sqrt{V(\Chi_B)} \geq 0.
$
This inequality implies that $W \geq 0$ and
$W^{-1}(0) \subset \underline{W}^{-1}(0)$. But $\underline{W}=0$ if and
only if $V(\Chi_B)=0$ (i.e., $\Chi_B=0$). Thus $W^{-1}(0)
\subset \{(\Chi_B,R_B): \Chi_B=0\}$. Conversely, on the set $\{(\Chi_B,R_B): \Chi_B=0\}$, $\Chi_B=0$ and hence $W=0$, and therefore $\{(\Chi_B,R_B): \Chi_B=0\} \subset W^{-1}(0)$. It follows that $W^{-1}(0)
= \{(\Chi_B,R_B): \Chi_B=0\}$ proving part (i).

For part (ii), note that for all $c>0$, $W_c \subset \{
\underline{W}(\Chi,R) \leq c\}$. Since the sublevel sets of $\underline{W}$ are compact and $R_B\in \sR_B$, a compact set, the set $W_c$ is bounded. Continuity of $W$ implies that $W_c$ is compact. 

For part (iii), it has already been shown that $W(\Chi_B,R_B) \geq \alpha^\star \sqrt{V(\Chi_B)}$. It also holds that $W=\sqrt{V(\Chi_B)} \left( \alpha + \sum_{i \in B}\hat\wdes_i^i\left(\mu_i(\Chi_B),R_i\right)\cdot e_1 \right)\leq \sqrt{V(\Chi_B)} \left( \alpha + \alpha \right)\leq 2\alpha \sqrt{V(\Chi_B)}$. \qed

\subsection{Proof of Lemma~\ref{lem:inequalities}}

We first compute inequalities for $\dot W_\tran$ and $\dot W_\rot$ for system~\eqref{eq:rel_coords_trans} and~\eqref{eq:rel_coords_rot}. We then combine them to derive \eqref{eq:ineq}.  
Consider unicycle $i \in B$. The dynamics of $\Chi_i$ in~\eqref{eq:rel_coords_trans} are split into two terms, for neighboring robots $j \in {\cal N}_i \cap A$ and $j \in {\cal N}_i \cap B$ respectively,
\begin{equation}\label{eq:split_dyn}
\begin{aligned}
\dot \Chi_i =&\sum_{j\in {\cal N}_i \cap A}a_{ij}\frac{(u_jR_je_1-u_iR_ie_1)}{\p}\\
&+\sum_{j\in {\cal N}_i \cap B}a_{ij}\frac{(u_jR_je_1-u_iR_ie_1)}{\p}.
\end{aligned}
\end{equation}
For simplicity of notation, we drop the arguments of $\hat\gdes_i(\Chi_i)$ and $\hat\gdes_i^i(\Chi_i,R_i)$. Adding and subtracting the term, 
\[
\frac{\sum_{j\in {\cal N}_i \cap B}a_{ij}(\hat\gdes_j-\hat\gdes_i)-\sum_{j\in {\cal N}_i \cap A}a_{ij}\hat\gdes_i}{\p}
\]
to~\eqref{eq:split_dyn} yields,
\[
\begin{aligned}
&\dot \Chi_i =\frac{\sum_{j\in {\cal N}_i \cap B}a_{ij}(\hat\gdes_j-\hat\gdes_i)-\sum_{j\in {\cal N}_i \cap A}a_{ij}\hat\gdes_i}{\p}\\
+&\frac{\sum_{j\in {\cal N}_i \cap B}a_{ij}(u_jR_je_1-u_iR_ie_1)}{\p}-\frac{\sum_{j\in {\cal N}_i \cap B}a_{ij}(\hat\gdes_j-\hat\gdes_i)}{\p}\\
+&\frac{\sum_{j\in {\cal N}_i \cap A}a_{ij}u_jR_je_1}{\p}+\frac{\sum_{j\in {\cal N}_i \cap A}a_{ij}(\hat\gdes_i-u_iR_ie_1)}{\p}\\
&=\frac{\sum_{j\in {\cal N}_i \cap B}a_{ij}(\hat\gdes_j-\hat\gdes_i)-\sum_{j\in {\cal N}_i \cap A}a_{ij}\hat\gdes_i}{\p}\\
+&\frac{\sum_{j\in {\cal N}_i \cap B}a_{ij}(u_jR_je_1-\hat\gdes_j)}{\p}-\frac{\sum_{j\in {\cal N}_i \cap B}a_{ij}(u_iR_ie_1-\hat\gdes_i)}{\p}\\
+&\frac{\sum_{j\in {\cal N}_i \cap A}a_{ij}u_jR_je_1}{\p}+\frac{\sum_{j\in {\cal N}_i \cap A}a_{ij}(\hat\gdes_i-u_iR_ie_1)}{\p}.
\end{aligned}
\]
Replacing $u_j$ and $u_i$ by the
assigned feedbacks in~\eqref{eq:CCP} and using the identity
$R_i \hat\gdes_i^i = \hat\gdes_i$ then,
\[
\dot \Chi_i=a_i(\Chi_B)+b_i(\Chi_B,R)+c_i(\Chi_B,R)+d_i(\Chi_A,R),
\]
where,
\[
\begin{aligned}
a_i(\Chi_B)&:=\frac{\sum_{j\in {\cal N}_i \cap B}a_{ij}(\hat\gdes_j-\hat\gdes_i)-\sum_{j\in {\cal N}_i \cap A}a_{ij}\hat\gdes_i}{\p}\\
b_i(\Chi_B,R)&:=\frac{\sum_{j\in {\cal N}_i \cap B}a_{ij}R_j((\hat\gdes_j^j \cdot e_1)e_1-\hat\gdes_j^j)}{\p}\\
&-\frac{\sum_{j\in {\cal N}_i \cap B}a_{ij}R_i((\hat\gdes_i^i \cdot e_1)e_1-\hat\gdes_i^i)}{\p}\\
c_i(\Chi_B,R)&:=\frac{\sum_{j\in {\cal N}_i \cap A}a_{ij}R_i(\hat\gdes_i^i-(\hat\gdes_i^i \cdot e_1)e_1)}{\p}\\
d_i(\Chi_A,R)&:=\frac{\sum_{j\in {\cal N}_i \cap A}a_{ij}(\hat\gdes_j^j \cdot e_1)R_je_1}{\p}.
\end{aligned} 
\]
The time derivative of $W_\tran=\sqrt{V(\Chi_B)}$ in~\eqref{eq:lyaps} yields, 
\begin{equation}\label{eqn:deriv}
\begin{aligned}
\dot W_\tran &= \frac{1}{2\sqrt{V(\Chi_B)}} \Bigg[\sum_{i\in B} \frac{\partial V(\Chi_B)}{\partial \Chi_i}(a_i(\Chi_B)+b_i(\Chi_B,R)\\
&+c_i(\Chi_B,R))\Bigg]+\frac{1}{2\sqrt{V(\Chi_B)}}\sum_{i\in B} \frac{\partial V(\Chi_B)}{\partial \Chi_i}d_i(\Chi_A,R).
\end{aligned}
\end{equation}
The derivative of the first term is considered in Claim~\ref{lem:sinteg}.
\begin{claim}\label{lem:sinteg} 
There exist gains $\gamma_i$ in~\eqref{eq:lyaps} and a negative definite function $\kdes(\Chi_B)$, homogeneous of degree three, such that $
\sum_{i\in B} \frac{\partial V(\Chi_B)}{\partial \Chi_i}a_i(\Chi_B) \leq \kdes(\Chi_B)$.  
\end{claim} 
The proof of Claim~\ref{lem:sinteg} is presented in Section~\ref{sec:claim} of this Appendix. Let the gains $\gamma_i$ be as in Claim~\ref{lem:sinteg}. The derivative of the remaining terms in the square brackets of~\eqref{eqn:deriv} satisfies,
\[
\begin{aligned}
&\sum_{i\in B} \frac{\partial V(\Chi_B)}{\partial \Chi_i}(b_i(\Chi_B,R)+c_i(\Chi_B,R))\\
&\leq \sum_{i\in B} \frac{1}{\p} \frac{\partial V(\Chi_B)}{\partial \Chi_i} \left[ \sum_{j\in {\cal N}_i \cap B}a_{ij}\left\|(\hat\gdes_j^j \cdot e_1)e_1-\hat\gdes_j^j\right\|\right.\\
&\left.+\sum_{j\in {\cal N}_i \cap B}a_{ij}\left\|(\hat\gdes_i^i \cdot e_1)e_1-\hat\gdes_i^i\right\|+\sum_{j\in {\cal N}_i \cap A}a_{ij}\left\|\hat\gdes_i^i-(\hat\gdes_i^i \cdot e_1)e_1\right\|\right]\\
&\leq \sum_{i\in B} \frac{1}{\p} \frac{\partial V(\Chi_B)}{\partial \Chi_i} \left[ \sum_{j\in {\cal N}_i \cap B}a_{ij}\left\|(\hat\gdes_j^j \cdot e_1)e_1-\hat\gdes_j^j\right\|\right.\\
&\left.+\sum_{j\in {\cal N}_i}a_{ij}\left\|(\hat\gdes_i^i \cdot e_1)e_1-\hat\gdes_i^i\right\| \right].
\end{aligned}
\]
We claim that $\|( \hat\gdes_i^i(\Chi_i,R_i) \cdot e_1 ) e_1 - \hat\gdes_i^i(\Chi_i,R_i) \| =
\left|\hat\gdes_i^i(\Chi_i,R_i) \cdot e_2 \right|$. Indeed, writing $\hat\gdes_i^i = (\hat\gdes_i^i
\cdot e_1) e_1 + \hat\gdes_i^i - (\hat\gdes_i^i \cdot e_1) e_1$, we have $ \hat\gdes_i^i
\cdot e_2 = ( \hat\gdes_i^i - (\hat\gdes_i^i \cdot e_1) e_1 ) \cdot e_2$.  Since
the vector $\hat\gdes_i^i - (\hat\gdes_i^i \cdot e_1) e_1$ is parallel to
$e_2$, $\left|( \hat\gdes_i^i - (\hat\gdes_i^i \cdot e_1) e_1 ) \cdot e_2 \right| =
\|\hat\gdes_i^i - (\hat\gdes_i^i \cdot e_1) e_1 \|$, so that $\left|\hat\gdes_i^i \cdot
e_2\right| = \|\hat\gdes_i^i - (\hat\gdes_i^i \cdot e_1) e_1 \|$. Then,
\[
\begin{aligned}
\sum_{i\in B} &\frac{\partial V(\Chi_B)}{\partial \Chi_i}(b_i(\Chi_B,R)+c_i(\Chi_B,R))\\
&\leq \sum_{i\in B} \frac{\bar a}{\p} \left\|\frac{\partial V(\Chi_B)}{\partial \Chi_i}\right\| \left( \sum_{j\in B}\left|\hat\gdes_j^j \cdot e_2 \right|+n\left|\hat\gdes_i^i \cdot e_2 \right|\right)
\end{aligned}
\]
where $\bar a=\max\{a_{ij}\}_{i,j \in \{1,\dots,n\}}$ which is homogeneous of degree three with respect to $\Chi_B$ since $\frac{\partial V(\Chi_B)}{\partial \Chi_i}$ is homogeneous of degree one and $\hat\gdes_i^i$ is homogeneous of degree two with respect to $\Chi_B$ for all $i \in B$. The last term in~\eqref{eqn:deriv} satisfies,
\begin{equation}\label{eq:last_term}
\begin{aligned}
&\frac{1}{2\sqrt{V(\Chi_B)}}\sum_{i\in B} \frac{\partial V(\Chi_B)}{\partial \Chi_i}d_i(\Chi_A,R)\\
&\leq \frac{1}{2\sqrt{V(\Chi_B)}} \sum_{i\in B} \frac{1}{\p}\frac{\partial V(\Chi_B)}{\partial \Chi_i}\sum_{j\in {\cal N}_i \cap A}a_{ij}(\hat\gdes_j^j \cdot e_1)R_je_1\\
&\leq \frac{1}{2\sqrt{V(\Chi_B)}} \sum_{i\in B} \frac{1}{\p} \left\|\frac{\partial V(\Chi_B)}{\partial \Chi_i}\right\|\sum_{j\in {\cal N}_i \cap A}a_{ij}\|\hat\gdes_j^j\|\\
&\leq \sum_{i\in B} \sup_{\Chi_B \in \sX_B} \left\{\frac{\bar a}{\p}\frac{1}{2\sqrt{V(\Chi_B)}} \right.\\
&\left.\left\|\frac{\partial V(\Chi_B)}{\partial \Chi_i}\right\| \right\}\sum_{j\in {\cal N}_i \cap A}\|\hat\gdes_j^j\|:=\Phi_{\tran}(\Chi_A,R).
\end{aligned}
\end{equation}
The bounded supremum of $\frac{1}{\sqrt{V(\Chi_B)}} \left\|\frac{\partial V(\Chi_B)}{\partial \Chi_i}\right\|$ exists because this term is homogeneous of degree 0 with respect to $\Chi_B$. Moreover, $\Chi_A=0$ implies that $\|\hat\gdes_j^j\|= 0$ for all $j \in A$ and hence $\Phi_{\tran}(0,R)=0$.
Everything together, \eqref{eqn:deriv} yields,
\begin{equation}\label{eqn:deriv2}
\begin{aligned}
\dot W_\tran &\leq \frac{1}{2\sqrt{V(\Chi_B)}}\left[\kdes(\Chi_B)+\sum_{i\in B} \frac{\bar a}{\p} \left\|\frac{\partial V(\Chi_B)}{\partial \Chi_i}\right\|\right.\\ 
& \left.\left( \sum_{j\in B}\left|\hat\gdes_j^j \cdot e_2 \right|+n\left|\hat\gdes_i^i \cdot e_2 \right|\right)\right]+\Phi_{\tran}(\Chi_A,R).
\end{aligned}
\end{equation}
Since $\kdes(\Chi_B)$ is homogeneous of degree three. We can write, 
\[
\begin{aligned}
\kdes(\Chi_B)&=\frac{\sqrt{V(\Chi_B)}V(\Chi_B)}{\sqrt{V(\Chi_B)}V(\Chi_B)}\kdes(\Chi_B)\\
&=\sqrt{V(\Chi_B)}V(\Chi_B)\kdes\left(\frac{\Chi_B}{\sqrt{V(\Chi_B)}}\right)\\
&=\sqrt{V(\Chi_B)}V(\Chi_B)\kdes\left(\mu(\Chi_B)\right).
\end{aligned}
\]
Analogous operations can be performed with the remaining term in the square bracket of~\eqref{eqn:deriv2} yielding,
\[
\begin{aligned}
&\dot W_\tran\leq \frac{V(\Chi_B)}{2}\Bigg[\kdes(\mu(\Chi_B)) +\sum_{i\in B} \frac{\bar a}{\p} \left\|\frac{\partial V(\mu(\Chi_B))}{\partial \Chi_i}\right\|\\
&\left. \left( \sum_{j\in B}\left|\hat\gdes_j^j(\mu_j(\Chi_B),R_j) \cdot e_2 \right|+n\left|\hat\gdes_i^i(\mu_i(\Chi_B),R_i) \cdot e_2 \right|\right)\right]\\
&+\Phi_{\tran}(\Chi_A,R).
\end{aligned}
\]
Since $\kdes$ is continuous and negative definite and $\mu(\Chi_B)$ lies on a compact set $S_1$, it follows that $\kdes(\mu(\Chi_B))/2$ has bounded maximum $-M_2<0$. Similarly, the function $\frac{\bar a}{\p} \left\|\frac{\partial V(\mu(\Chi_B))}{\partial \Chi_i} \right\|$ has a maximum. Letting $M_1:=n \max_{\stackrel{\prj \in S_1 }{i \in B}} \frac{\bar a}{\p} \left\|\frac{\partial V(\prj)}{\partial \Chi_i}\right\|$ yields,
\begin{equation}\label{eq:Wtrandot}
\begin{aligned}
\dot W_\tran \leq &V(\Chi_B)\left[-M_2+\frac{M_1}{2n}\sum_{i\in B}\left( \sum_{j\in B}\left|\hat\gdes_j^j(\mu_j(\Chi_B),R_j) \cdot e_2 \right|\right.\right.\\
&\left.\left.+n\left|\hat\gdes_i^i(\mu_i(\Chi_B),R_i) \cdot e_2 \right|\right)\right]+\Phi_{\tran}(\Chi_A,\sR)\\
\leq &V(\Chi_B)\left[-M_2+\frac{M_1}{2n}\sum_{i\in B}\left( n\left|\hat\gdes_i^i(\mu_i(\Chi_B),R_i) \cdot e_2 \right|\right.\right.\\
&\left.\left. +n\left|\hat\gdes_i^i(\mu_i(\Chi_B),R_i) \cdot e_2 \right|\right)\right]+\Phi_{\tran}(\Chi_A,\sR)\\
\leq &V(\Chi_B)\left[-M_2+M_1\sum_{i\in B}\left|\hat\gdes_i^i(\mu_i(\Chi_B),R_i) \cdot e_2 \right|\right]\\
&+\Phi_{\tran}(\Chi_A,\sR).\\
\end{aligned}
\end{equation}
This proves the first inequality. We now turn to the second. Recall the definition of $W_\rot$,
$
W_\rot(\Chi_B,R_B) = \sum_{i \in B} \hat\wdes_i^i(\Chi_i,R_i) \cdot e_1.
$
The time derivative of $W_\rot$ along the vector field
in~\eqref{eq:rel_coords_trans}-\eqref{eq:rel_coords_rot} is
$
\dot W_\rot = \sum_{i \in B} \left( \frac{d}{dt} \hat\wdes_i^i \right)
  \cdot e_1.
$
To express $(d / dt) \hat\wdes_i^i$, recall that $\hat\wdes_i^i(\Chi_i,R_i) =
R_i^{-1} \hat \wdes_i(\Chi_i)$. Then,
$
\frac d {dt} \hat\wdes_i^i = \left( \frac d {dt} R_i^{-1} \right)
\hat\wdes_i + R_i^{-1} \frac {d \hat\wdes_i}{dt}.
$
We will denote the derivative of $\hat \wdes_i(\Chi_i)=\p \Chi_i$ by,
\[
\begin{aligned}
\hdes_i&(\Chi,R) :=(d/dt)\hat\wdes_i(\Chi_i)\\
&=\p\left(a_i(\Chi_B)+b_i(\Chi_B,R)+c_i(\Chi_B,R)+ d_i(\Chi_A,R)\right)
\end{aligned}
\]
where the first three terms are homogeneous of degree two with respect to $\Chi_B$ and the last term is homogeneous of degree two with respect to $\Chi_A$. Consistently with our notational convention, we will let $\hdes_i^i(\Chi,R):=R_i^{-1}
\hdes_i(\Chi,R)$. Returning to the derivative of $\hat\wdes_i^i$, we have
\[
\begin{aligned}
\frac d {dt} \hat\wdes_i^i &= - (\omega_i)^\times R_i^{-1} \hat\wdes_i(\Chi_i)
+ R_i^{-1} \hdes_i(\Chi,R) \\
&= -\begin{bmatrix}
  0  & -\omega_i \\
\omega_i & 0 
\end{bmatrix} \hat\wdes_i^i(\Chi_i,R_i) + \hdes_i^i(\Chi,R).
\end{aligned}
\]
We substitute the above identity in the expression for $\dot W_\rot$,
\[
\begin{aligned}
\dot W_\rot &= \sum_{i\in B} \left(-e_1^\top \begin{bmatrix}
  0  & -\omega_i \\
\omega_i & 0 
\end{bmatrix} \hat\wdes_i^i(\Chi_i,R_i) + \hdes_i^i(\Chi,R)\cdot e_1 \right) \\
&=\sum_{i \in B} \left( (\hat\wdes_i^i(\Chi_i,R_i)\cdot e_2)\omega_i + \hdes_i^i(\Chi,R)\cdot e_1 \right).\\
\end{aligned}
\]
Substituting the feedback $\omega_i=-k_1
(\hat \wdes_i^i(\Chi_i,R_i) \cdot e_2)$ and taking norms, we arrive at the inequality
\[
\dot W_\rot  \leq \sum_{i\in B} \big[ - k_1 \left| \hat\wdes_i^i(\Chi_i,R_i) \cdot e_2\right|^2
  + \hdes_i^i(\Chi,R) \cdot e_1 \big].
\]
This gives,
\[
\dot W_\rot\leq \left[ -k_1 \sum_{i \in B} \left| \hat\wdes_i^i(\Chi_i,R_i) \cdot e_2\right|^2 + \ell(\Chi_B,R)\right]+\Phi_{\rot}(\Chi_A,R)\\
\]
where
\[
\ell(\Chi_B,R):=\sum_{i \in B} \p R_i^\top\left(a_i(\Chi_B)+b_i(\Chi_B,R)+c_i(\Chi_B,R)\right) \cdot e_1
\]
and $\Phi_{\rot}(\Chi_A,R):=\sum_{i \in B} \p R_i^\top d_i(\Chi_A,R) \cdot e_1$. Note that $\sum_{i \in B} \left| \hat\wdes_i^i(\Chi_i,R_i) \cdot e_2\right|^2$ and $\ell(\Chi_B,R)$ are homogeneous of degree two with respect to $\Chi_B$. The function $\Phi_{\rot}(\Chi_A,R)$ does not depend on $\Chi_B$ and $\Phi_{\rot}(0,R)=0$. This yields,
\[
\begin{aligned}
\dot W_\rot \leq& V(\Chi_B)\left[ - k_1 \sum_{i\in B}\left| \hat\wdes_i^i(\Chi_i/\sqrt{V(\Chi_B)},R_i) \cdot e_2\right|^2 \right.\\
&\left. + \ell(\Chi_B/\sqrt{V(\Chi_B)},R) \right]+\Phi_{\rot}(\Chi_A,R)\\
&\leq V(\Chi_B) \left[ - k_1 \sum_{i\in B}\left| \hat\wdes_i^i(\mu_i(\Chi_B),R_i) \cdot e_2\right|^2\right.\\
&\left.+ \ell(\mu(\Chi_B),R) \right]+\Phi_{\rot}(\Chi_A,R).
\end{aligned}
\]
$\left| \ell(\mu(\Chi_B),R) \right|$ has a bounded supremum. Letting $M_3 = \sup_{(\prj,R)\in S_1 \times \sR} \left(\left| \ell(\prj,R) \right| \right)$, we conclude that,
\begin{equation}\label{eq:Wrotdot}
\begin{aligned}
\dot W_\rot \leq& V(\Chi_B)\left[ -k_1 \sum_{i\in B}\left|
  \hat\wdes_i^i(\mu_i(\Chi_B),R_i) \cdot e_2\right|^2 +
M_3 \right]\\
&+\Phi_{\rot}(\Chi_A,R).
\end{aligned}
\end{equation}

By using the inequalities \eqref{eq:Wtrandot} and \eqref{eq:Wrotdot} we now bound the derivative of $W$ to derive \eqref{eq:ineq}. Notice that 
\[
\begin{aligned}
\dot W&=\alpha \dot W_\tran+\dot W_\rot \\
\leq& V(\Chi_B) \left[-\alpha M_2 + \alpha M_1 \sum_{i \in B} \left|
\hat\gdes_i^i(\mu_i(\Chi_B),R_i) \cdot e_2 \right|\right.\\
&\left. -k_1 \sum_{i \in B} \left|
  \hat\wdes_i^i(\mu_i(\Chi_B),R_i) \cdot e_2\right|^2 +
M_3 \right]+\Phi(\Chi_A,R),
\end{aligned}
\]
where $\Phi(\Chi_A,R):=\alpha \Phi_{\tran}(\Chi_A,R)+\Phi_{\rot}(\Chi_A,R)$. 

Choose $\alpha>3M_3/M_2$. This implies,
\begin{equation}\label{eq:der_W}
\begin{aligned}
\dot W &\leq V(\Chi_B) \left[-2 M_3 + \alpha M_1 \sum_{i\in B} \left|
\hat\gdes_i^i(\mu_i(\Chi_B),R_i) \cdot e_2 \right|\right.\\
&\left.-k_1 \sum_{i\in B} \left|
  \hat\wdes_i^i(\mu_i(\Chi_B),R_i) \cdot e_2\right|^2 \right]+\Phi(\Chi_A,R).
\end{aligned}
\end{equation}
Since $\hat\wdes_i^i(\Chi_i,R_i)$ is homogeneous with respect to $\Chi_i$, we have, $\hat\wdes_i^i(\mu_i(\Chi_B),R_i)=\frac{\sqrt{\|\hat\gdes_i^i(\mu_i(\Chi_B),R_i)\|}}{\|\hat\gdes_i^i(\mu_i(\Chi_B),R_i)\|}\hat\gdes_i^i(\mu_i(\Chi_B),R_i)$. Plugging the last expression into~\eqref{eq:der_W} yields
\[
\begin{aligned}
\dot W \leq& V(\Chi_B) \left[-2 M_3 + \alpha M_1 \sum_{i\in B} \left|
\hat\gdes_i^i(\mu_i(\Chi_B),R_i) \cdot e_2 \right| \right.\\
&\left.-k_1 \sum_{i\in B}
  \left(\frac{\sqrt{\|\hat\gdes_i^i(\mu_i(\Chi_B),R_i)\|}}{\|\hat\gdes_i^i(\mu_i(\Chi_B),R_i)\|}\hat\gdes_i^i(\mu_i(\Chi_B),R_i) \cdot e_2\right)^2 \right]\\
 &+\Phi(\Chi_A,R)\\  
\leq& V(\Chi_B) \left[-2 M_3 + \alpha M_1 \sum_{i\in B} \left|
\hat\gdes_i^i(\mu_i(\Chi_B),R_i) \cdot e_2 \right|\right.\\
&\left.-k_1 \sum_{i\in B} \frac{1}{\|\hat\gdes_i^i(\mu_i(\Chi_B),R_i)\|}\left|\hat\gdes_i^i(\mu_i(\Chi_B),R_i) \cdot e_2\right|^2 \right]\\
 &+\Phi(\Chi_A,R).
\end{aligned}
\]
Since $\hat\gdes_i^i(\mu_i(\Chi_B),R_i)$ is a continuous function of its arguments and $\mu_i(\Chi_B)$ is compact, $\left|\hat\gdes_i^i(\mu_i(\Chi_B),R_i)\right|$ has a maximum $M_4$. This implies,
\[
\begin{aligned}
\dot W &\leq V(\Chi_B) \left[-2 M_2 + \alpha M_1 \sum_{i\in B} \left|
\hat\gdes_i^i(\mu_i(\Chi_B),R_i) \cdot e_2 \right| \right.\\
&\left.-k_1 \sum_{i\in B}
  \frac{1}{M_4}\left|\hat\gdes_i^i(\mu_i(\Chi_B),R_i) \cdot e_2\right|^2 \right]+\Phi(\Chi_A,R).
\end{aligned}
\]
Denote $\bb_i(\mu_i(\Chi_B),R_i):=\left|\hat\gdes_i^i(\mu_i(\Chi_B),R_i) \cdot e_2\right|$, and $\bb:=(\bb_i(\mu_i(\Chi_B),R_i))_{i\in B}$. Then,
\[
\begin{aligned}
\dot W &\leq V(\Chi_B) \left[-2 M_2 + \alpha M_1 \one\trans \bb-\frac{k_1}{M_4}\left|\bb\right|^2 \right]+\Phi(\Chi_A,R)\\
&= V(\Chi_B) \big[ \one\trans \ \ \bb\trans\big] 
\begin{bmatrix}
  \frac{-2 M_2} n I  & \alpha \frac{M_1}{2} I \\
\alpha \frac{M_1}{2} I & \frac{-\kq}{M_4} I 
\end{bmatrix}
\begin{bmatrix} 
\one \\ \bb
\end{bmatrix}+\Phi(\Chi_A,R).
\end{aligned}
\]
There exists $\kq^\star>0$ such that choosing $\kq>\kq^\star$, the matrix above is negative definite and therefore the first term satisfies,
\begin{equation}\label{eq:Lyap_ineq}
V(\Chi_B) \big[ \one\trans \ \ \bb\trans\big] 
\begin{bmatrix}
  \frac{-2 M_2} n I  & \alpha \frac{M_1}{2} I \\
\alpha \frac{M_1}{2} I & \frac{-\kq}{M_4} I 
\end{bmatrix}
\begin{bmatrix} 
\one \\ \bb
\end{bmatrix} \leq -\sigma V(\Chi_B),
\end{equation}
$\sigma>0$. This concludes the proof of Lemma~\ref{lem:inequalities}. \qed

\subsection{Proof of Claim~\ref{lem:sinteg}}\label{sec:claim}
Recalling that $V(\Chi_B)=\gamma_i\Chi_i^\top \Chi_i$ with $\Chi_i=\hat \wdes_i/\p$ and defining $b_{ij}:=\frac{a_{ij}}{\p^2}$, it holds that, 
\[
\begin{aligned}
\sum_{i\in B} &\frac{\partial V(\Chi_B)}{\partial \Chi_i}a_i(\Chi_B)=2\sum_{i\in B} \gamma_i \frac{\hat \wdes_i}{\p} \cdot a_i(\Chi_B)\\ 
\leq&2\sum_{i\in B}\gamma_i \hat\wdes_i \hspace{-0.5mm}\cdot \hspace{-1mm}
\left(\sum_{j \in {\cal N}_i \cap B} \hspace{-1mm}b_{ij} (\|\hat \wdes_j\|\hat \wdes_j-\|\hat \wdes_i\|\hat \wdes_i)-\hspace{-2.4mm}\sum_{j \in {\cal N}_i \cap A} b_{ij}\|\hat \wdes_i\|\hat \wdes_i \right) \\
\leq&2\sum_{i\in B}\gamma_i \hspace{-1mm} \left(\sum_{j \in {\cal N}_i \cap B}\hspace{-1mm} b_{ij} (-\|\hat \wdes_i\|^3+ \|\hat \wdes_j\|\hat \wdes_j \cdot \hat \wdes_i)-\hspace{-3.5mm}\sum_{j \in {\cal N}_i \cap A} \hspace{-1mm}b_{ij} \|\hat \wdes_i\|^3 \right)\\
\leq& \sum_{i\in B}\gamma_i\sum_{j \in {\cal N}_i \cap B} b_{ij}\left(-\frac{4}{3}\|\hat \wdes_i\|^3+ \frac{4}{3} \|\hat \wdes_j\|^3\right)\\
&+\sum_{i\in B} \gamma_i\sum_{j \in {\cal N}_i \cap B} b_{ij} \left(-\frac{2}{3}\|\hat \wdes_i\|^3+ 2\|\hat \wdes_j\|\hat \wdes_j \cdot \hat \wdes_i-\frac{4}{3}\|\hat \wdes_j\|^3\right)\\
&-2\sum_{i\in B}\gamma_i\sum_{j \in {\cal N}_i \cap A} b_{ij} \|\hat \wdes_i\|^3.
\end{aligned}
\] 
The first term equals $\frac{4}{3} \gamma^\top M \bar h$ with $\bar h:=(\|\hat \wdes_i\|^3)_{i \in B}$. $M$ is the $(r \times r)$-matrix whose $(i,j)$-th component is $\sum_{k\in {\cal N}_{i}\cap B}b_{ik}$ for $i=j$, $b_{ij}$ for $j\in {\cal N}_{i} \cap B$ and zero otherwise for $i,j \in \{1,\dots,r\}$ where it is assumed without loss of generality that $B=\{1,\dots,r\}$. Choose $\gamma=(\gamma_1,\dots,\gamma_n)$ as the left eigenvector associated to the zero eigenvalue of $M$. Since $B$ corresponds to a collection of strongly connected components with no links from one to the other, the zero eigenvalue is unique and all components of $\gamma$ are positive (see Proposition D.5 in~\cite{hatanaka2015}). Therefore, 
\[
\begin{aligned}
&\sum_{i\in B} \frac{\partial V(\Chi_B)}{\partial \Chi_i}a_i(\Chi_B)\\
\leq& \sum_{i \in B} \gamma_i \sum_{j \in {\cal N}_i \cap B} b_{ij} \left(-\frac{2}{3}\|\hat \wdes_i\|^3+ 2\|\hat \wdes_j\|\hat \wdes_j \cdot \hat \wdes_i-\frac{4}{3}\|\hat \wdes_j\|^3\right)\\
&-2\sum_{i \in B} \gamma_i \sum_{j \in {\cal N}_i \cap A} b_{ij} \|\hat \wdes_i\|^3=:\kdes(\Chi_B).
\end{aligned}
\]
The term 
\[
\begin{aligned}
\kdes_1(\Chi_B)&:=\sum_{i \in B} \gamma_i \hspace{-2mm}\sum_{j \in {\cal N}_i \cap B} \hspace{-2mm}b_{ij} \hspace{-1mm}\left(\hspace{-1mm}-\frac{2}{3}\|\hat \wdes_i\|^3+2\|\hat \wdes_j\|\hat \wdes_j \cdot \hat \wdes_i-\frac{4}{3}\|\hat \wdes_j\|^3\right)\\
&\leq \sum_{i \in B} \gamma_i \hspace{-2mm}\sum_{j \in {\cal N}_i \cap B} \hspace{-2mm}b_{ij} \hspace{-1mm}\left(-\frac{2}{3}\|\hat \wdes_i\|^3+2\|\hat \wdes_i\|\|\hat \wdes_j\|^2-\frac{4}{3}\|\hat \wdes_j\|^3\right)
\end{aligned}
\]
is less than or equal to zero with equality only when $\hat\wdes_i= \hat\wdes_j$ for all $i,j \in B$ and as such $\kdes(\Chi_B)$ is less than or equal to zero with equality only when $\hat\wdes_i= \hat\wdes_j$ for all $i,j \in B$. 

Now we prove that $\kdes(\Chi_B) = 0$ only if $\hat \wdes_i=0$ for all robots $i \in B$. In the case that $A$ is not empty, the inequality $\kdes(\Chi_B) \leq -2\sum_{i \in B}\gamma_i\sum_{j \in {\cal N}_i \cap A} b_{ij} \|\hat \wdes_i\|^3$ implies $\kdes(\Chi_B)=0$ only if $\hat \wdes_i=0$ for any $i \in B$ with a neighbor in $A$. As such, by the previous arguments, $\kdes(\Chi_B)=0$ only if $\hat \wdes_i=0$ for all $i \in B$. On the other hand, if $A$ is empty, then $B$ is isolated and strongly connected. Therefore $\kdes(\Chi_B)= \kdes_1(\Chi_B)$ is equal to zero only if $\kdes_1(\Chi_B)=0$ which is the case only if $\hat \wdes_i=\hat \wdes_j$ for all $i,j \in B$. This implies that $(L \otimes I_2) x\in \Span\{\one \otimes e_1,\one \otimes e_2\}$. Since $B$ is a strongly connected component there exists a unique vector $\bar \gamma$ (with positive entries) such that $\bar \gamma^\top (L \otimes I_2) = 0$. Since $\bar \gamma^\top (L \otimes I_2)x = \bar \gamma^\top \one \otimes (\alpha e_1+\beta e_2)$ for some $\alpha, \beta \in \Re$, it holds that $\bar \gamma^\top \one \otimes (\alpha e_1+\beta e_2)=0$. Since all entries of $\bar \gamma$ are positive, this implies $\alpha=\beta=0$ and $(L \otimes I_2)x=0$. Therefore $x \in \Span\{\one \otimes e_1,\one \otimes e_2\}$ or, equivalently, that $\hat \wdes_i=0$ for all $i\in B$. 

Therefore $\kdes(\Chi_B)=0$ only if $\Chi_i=0$ for all $i\in B$ and as such $\kdes(\Chi_B)$ is negative definite. Note that $\kdes(\Chi_B)$ is homogeneous of degree three with respect to $\Chi_B$ because $\hat \wdes_i$ is homogeneous of degree one with respect to $\Chi_B$ for all $i \in B$. This completes the proof of the claim.

\end{document}